\documentclass[11pt,a4paper,sunlimits]{article}

\usepackage{indentfirst}
\usepackage{amsmath}
\usepackage{amsthm}
\usepackage{amsfonts}
\usepackage{amssymb}
\usepackage{amscd}
\usepackage{textcomp}
\usepackage{graphicx}

\usepackage[all]{xy}
\newenvironment{dem}{\noindent{\bf Demonstration}~:}%
{\unskip\hfill\null\nobreak\hfill\carre\vskip1em\par}
\newcommand{\carre}{\rule{1ex}{1ex}}

\newtheorem{Theorem}{Theorem}[section]
\newtheorem{lemma}[Theorem]{Lemma}

\newtheorem{cor}[Theorem]{Corollary}
\newtheorem{Definition}[Theorem]{Definition}

 \def\fd {\hspace{0.35cm} \raise
-0.5mm\hbox{$\blacksquare$}\\}

\def\qed{{\null\hfill\ \raise3pt\hbox{\framebox[0.1in]{}}}\break\null}

\def \g{\mathfrak}

\def\a{\alpha}
\def\b{\beta}

\def\d{\delta}
\def\e{\varepsilon}

\def\ga{\gamma}
\def\G{\Gamma}

\def\l{\lambda}
\def\o{\omega}
\def\O{\Omega}

\def\s{\sigma}
\def\S{\Sigma}
\def\T{\Theta}
\def\t{\theta}

\def\cB{{\cal B}}

\def\cD{{\cal D}}

\def\cF{{\cal F}}

\def\cM{{\cal M}}
\def\cN{{\cal N}}
\def\cO{{\cal O}}

\def\cR{{\cal R}}
\def\cS{{\cal S}}

\def\cU{{\cal U}}
\def\cV{{\cal V}}
\def\cW{{\cal W}}

\def\N{\mathbb N}
\def\Z{\mathbb Z}

\def\R{\mathbb R}
\def\C{\mathbb C}

\def\me{\medskip}
\def\no{\noindent}
\def\dis{\displaystyle}

\begin{document}

\author{Pascale Harinck\footnote{Ecole Polytechnique, CMLS - CNRS UMR 7640, Route de Saclay 91128 Palaiseau C\' edex, harinck@math.polytechnique.fr  }}
\title{Regularity of some invariant  distributions on nice symmetric pairs}
\date{}
\maketitle

\centerline{\bf Abstract}
\vskip 0,5cm

J.~Sekiguchi determined the semisimple symmetric pairs $(\g g, \g h)$, called nice symmetric pairs,  on which there is no non-zero invariant eigendistribution with singular support. On such pairs, we study regularity of invariant distributions annihilated by a polynomial of the Casimir operator. We deduce that invariant eigendistributions on ($\mathfrak{gl}(4,\R),\mathfrak{gl}(2,\R)\times \mathfrak{gl}(2,\R)$)  are   locally integrable functions.\me

\medskip


\section*{Introduction}
Let  $G$ be a reductive  group such that $\textrm{Ad}(G)$ is connected. Let   $\s$  be  an involutive automorphism of $G$. We denote by the same letter  $\s$ the corresponding involution on the Lie algebra $\g g$ of $G$. Let     $\g g=\g h\oplus \g q$ be the decomposition into $+1$ and $-1$ eigenspaces with respect to $\s$. Then $(\g g, \g h)$ is called a reductive symmetric pair  (or semisimple when $\g g$ is semisimple). Let    $H$ be the group of fixed points of $\s$ in $G$.

  In \cite{Se}, J. Sekiguchi  describes semisimple symmetric pairs  on which there is no non-zero invariant eigendistribution with support in  $\g q-\g q^{reg}$ where $\g q^{reg}$ is the set of semisimple regular elements of $\g q$. These pairs, called nice symmetric pairs, are characterized by a property on distinguished nilpotent elements and we can generalize this notion to reductive pairs (Definition \ref{defdistingue}). Our main result is the following . Let  $\o$ be the Casimir polynomial of $\g q$ and $\partial(\o)$ the corresponding differential operator on  $\g q$.

\begin{Theorem} Let $(\g g, \g h)$ be a nice reductive symmetric pair. Let $\cV$ be an $H$- invariant open subset of $\g q$. Let $\T$ be an $H$-invariant distribution on  $\cV$ such that
\begin{enumerate}\item There exists  $P\in\C[X]$ such that $P\big(\partial(\o)\big) \T=0$,
\item There exists $F\in L^1_{loc}(\cV)^H$ such that  $\T=F$ on  $\cV\cap\g q^{reg}$.
\end{enumerate}
Then $\T=F$ as distribution on  $\cV$. \end{Theorem}

 In \cite{GL}, E. Galina and Y. Laurent obtained  stronger results on invariant distributions on nice symmetric pairs  by different methods based on algebraic properties of $\cD$-modules.  They proved that any invariant distribution on a nice pair which is annihilated by a finite codimensional ideal of the algebra of $H$-invariant differential operators with constant coefficients on $\g q$ is a locally integrable function (\cite{GL} Corollary 1.7.6).\me

Our approach uses  properties of distributions. Assuming that  $S=\Theta-F  $ is non-zero, we are led to a contradiction. By the work of G. van Dijk (\cite{vd}) and J. Sekiguchi (\cite{Se} ), we can  adapt  the descent method of Harish-Chandra. Thus, we construct a non-zero distribution $\tilde{S}$ defined on a neighborhood $W$ of $0$ in  $\R^r\times\R^m$ with support in $(\{0\}\times \R^m)\cap W$ such that there exist a locally integrable function $\tilde{F}$ on  $W$ and a differential operator $D$, which is obtained from radial parts of  $\partial(\o)$ near semisimple elements and nilpotent elements,  satisfying   $P(D)\tilde{S}=P(D)\tilde{F}$. Using the method developed by M. Atiyah in   \cite{A}, one studies the   degree of singularity along $\{0\}\times \R^m$ of different distributions in this equation. One   deduces that $\tilde{S}=0$ and thus   a contradiction.\me

In the last section, we  complete the results of  \cite{HJ}  on the nice symmetric pair  $(\g gl(4,\R), \g gl(2,\R)\times\g gl(2,\R))$ and deduce that any invariant eigensdistribution for a regular character on this  pair is given by a locally integrable function.\me

\section{Notation}

 Let  $M$  be a smooth variety. Let  $C^{\infty}(M)$ be the space of smooth functions on $M$,
${\cal D}(M)$   the subspace of compactly supported smooth functions, $L^1_{loc}(M)$  the space of locally integrable functions on $M$, endowed with their standard topology and 
${\cal D}'(M)$  the space of distributions on  $M$.

For a group $G$  acting on $M$, one denotes by   $\cF^G$  the points of $\cF$ fixed by $G$ for each  space $\cF$ defined as above.

If $N\subset M$ and if  $f$ is a function defined on $M$, one denotes by $f _{/N}$
its restriction to $N$. 

If $V$  is a finite dimensional real vector space then $V^{*}$ is its algebraic dual 
and $V_{\Bbb C}$ is its complexified vector space.The symmetric algebra $S[V]$ of $V$
can be identified  to the space  $\R[V^*]$ of polynomial functions on $V^*$ with real coefficients  and to the space of  differential operators with real constant coefficients on $V$.  Similary, one has   $S[V_\C]=\C[V^*]$  and this algebra  can be identified  to the space   of  differential operators with complex constant coefficients on   $V_\C$. If  $u\in S[V]$ (resp. 
$S[V_\C]$), then  $\partial(u)$ will denote the corresponding differential operator.\me

Let  $G$ be a  reductive group such that $\textrm{Ad}(G)$ is connected, and  $\s$   an involution on $G$. This defines an involution, denoted by the same letter $\s$ on the Lie algebra $\g g$ of $G$.    Let  $\g g=\g h\oplus \g q$ be the direct decomposition of $\g g$ into the  $+1$ and $-1$ eigenspaces of $\s$. Then      $(\g g, \g h)$ is called a reductive symmetric pair. Let    $H$ be the   subgroup of  fixed points of $\s$ in $G$.\me

Let  $\g c_{\g g}$ be the center of $\g g$ and  $\g g_s$   its derived algebra. We set  $$\g c_{\g q}=\g c_{\g g}\cap \g q\textrm{ and } \g q_s=\g g_s\cap \g q.$$

If  $x$ is an element  of $\g g$ and $\g r$ is a subspace  of $\g g$, we denote by $\g r_x$ the centralizer of $x$ in $\g r$.\me

We fix a non-degenerate bilinear form $B$ on $\g g$ which is equal to the Killing form on $\g g_s$. Then $\o(X)=B(X,X)$ is the Casimir polynomial of $\g q$. 

\section{ Transfer of distributions and differential operators}\label{descente}

We recall results of  (\cite{vd} sections 2 and 3) and (\cite{Se} section (3.2)) on restriction of distributions and radial parts of differential operators. Their proofs are similar to (\cite{hc} or \cite{va} Part I, chapter 2 ).\me

Let $x_0\in \g q_s$. Let  $U$ be a linear subspace of $ \g q$ such that $\g q=U\oplus [x_0,\g h]$ and $V$  be a linear subspace of $\g h$ such that $\g h=V\oplus \g h_{x_0}$. Consider the open subset $^{\backprime}U=\{Z\in U; U+ [x_0+Z,\g h]=\g q\}$ containing $0$. 
Then the map 
$\Psi$ from $  H\times ^{\backprime}U $ to $ \g q$ defined by $\Psi(h,u)=h\cdot(x_0+u)$ 
is a submersion. In particular,  $\O=\Psi(H\times^{\backprime}U)$ is an open $H$-invariant subset of $\g q$ containing $x_0$. We fix an Haar measure  $dh$ on $H$ and we  denote by $du$  (respectively $dx$) the Lebesgue measure on  $U$ (respectively $\g q$). The submersion $\Psi$ induces a continuous surjective map  $\Psi_{\star}$ from $\cD(H\times^{\backprime}U)$ onto $\cD(\O)$ such that, for any  $F\in L_{loc}^1(\g q)$ and any $f\in\cD(H\times  ^{\backprime}U)$, one has
  $$\int_{H\times U} F\circ \Psi(h,u) f(h,u) dh\;du=\int_{\g q} F(x) \Psi_{\star}(f)(x) dx.$$

\begin{Theorem}\label{rest} For $T\in\cD'(\O)^H$ there exists a unique distribution $\cR es_U T$ defined on $^{\backprime}U$, called the restriction of $T$ to $^{\backprime}U$ with respect to $\Psi$, such that for any $f\in \cD(H\times  ^{\backprime}U)$, one has
$$<T,\Psi_{\star}(f)>=<\cR es_U T,p_{\star}(f)>$$
where $\dis p_{\star}(f)\in\cD(U)$ is defined by $\dis p_{\star}(f)(u)=\int_H f(h,u) dh$.

This restriction satisfies the following properties:
\begin{enumerate}
\item If $U$ is stable under the action of a subgroup $H_0$ of $H$ then  $\cR es_U T$ is $H_0$-invariant.
\item $x_0+\textrm{ supp }(\cR es_U T)\subset \textrm{ supp }(T)\cap(x_0+^{\backprime}U)$.

\item If  $F\in L^1_{loc}(\O)^H$ then  $\cR es_U F$ is the locally integrable function on   $^{\backprime}U$ defined by $\cR es_U F(u)=F(x_0+u)$ .
\item If $\cR es_U T=0$ then $T=0$ on   $\O$.
\end{enumerate}
\end{Theorem}

\begin{Theorem}\label{rad} Let  $D$ be a $H$-invariant differential operator on $\g q$. Then there exists a differential operator  $\cR ad_U(D)$, called the radial part of $D$ with respect to $\Psi$,  defined on  $^{\backprime} U$ such that for any $f\in \cD(\O)^H$, one has  $(D\cdot f)(x_0+u)=\cR ad_U(D)\cdot \cR es_U f(u)$ for  $u\in ^{\backprime}U$.

Morever, for any $T\in\cD'(\O)^H$, one has
$$\cR es_U(D\cdot T)=\cR ad_U(D)\cdot \cR es_U(T).$$
\end{Theorem} 

\section{Semisimple elements}\label{semisimple}

We recall that a Cartan subspace of  $\g q$ is a maximal abelian subspace of $\g q$  consisting of semisimple elements.

 If  $\g r=\g q$ or $\g q_s$, we denote by  $\cS(\g r)$ the set of semisimple elements of  $\g r$.\me

Let $\g a$ be a  Cartan subspace of  $\g q$. If  $\l\in\g g_\C^*$, we set  $$\g g_\C^\l=\{X\in\g g_\C; [A,X]=\l(A) X \textrm{ for any } A\in\g a_\C\}$$ and $$\S(\g a)=\{\l\in \g g_\C^*; \g g_\C^\l\neq\{0\}\}.$$
Then $\S(\g a)$ is the root system of ($\g g_\C, \g a_\C)$.\\
An element  $X$ of $\cS(\g q)$ is $\g q$-regular (or regular) if its centralizer $\g q_X$ in  $\g q$   is a Cartan subspace. If $X\in \g a$ then $X$ is regular if and only if  $\l(X)\neq 0$ for all $\l\in\S(\g a)$.
We denote by  $\g q^{reg}$ the open dense subset of semisimple regular elements of $\g q$. \me

Let $A_0\in\cS(\g q )$. Its centralizer $\g z=\g g_{A_0}$ in $\g g$ is a reductive   $\s$-stable Lie subalgebra of $\g g$. We denote by $\g c$ its center and by  $\g z_s$ its derived algebra. We set 
$$\g c^-=\g c\cap\g q, \quad\g c^+=\g c\cap\g h, \quad\g z_s^-=\g z_s\cap\g q \quad \textrm{and}  \quad\g z_s^+=\g z_s\cap\g h.$$ 
The pair  $(\g z_s,\g z_s^+)$ is a semisimple symmetric subpair of  $(\g g_s,\g h_s)$ which is  equal to $(\g g_s,\g h_s)$ if $A_0\in \g c_{\g q}$ . Let $H_s^+$ be the analytic subgroup of $H$ with Lie algebra $\g z_s^+$.\me

We assume that $A_0\notin\g c_{\g q}$. We take a Cartan subspace  $\g a$ of $\g q$ containing   $A_0$   and consider  the corresponding  root system $\S=\S(\g a)$. We fix  a positive system  $\S^+$ of $\S$. For any $\l\in\S^+$, we choose a $\C$-basis  $X_{\l,1},\ldots X_{\l,m_\l}$ of $\g g_\C^\l$ such that $B(X_{\l,i},\s(X_{\l,j}))=-\delta_{i,j}$ for $i,j\in\{1,\ldots, m_\l\}$. Let $\S_1^+=\{\l\in\S^+; \l(A_0)\neq 0\}$. We set  
$$ V_\C^\pm=\sum_{\l\in\S_1^+} \sum_{j=1}^{m_\l} \big(X_{\l,j}\pm\s(X_{\l,j})\big), \quad V^+=V^+_\C\cap\g h, \quad V^-=V^-_\C\cap \g q.$$
We have the decompositions  $\g h=\g z^+ \oplus V^+$ and $\g q=\g z^-\oplus V^-$, with $dim\; V^+=dim\; V^-$ and  $[A_0,\g h]=V^-$.

If $Z_0\in\g z^-$, we define the map $ \eta_{Z_0}$ from $ V^+\times \g z^-$ to $ \g q$ by $ \eta_{Z_0}(v, Z) = Z+[v,A_0+Z_0]$.
Then  $\eta_{ 0}$ is a bijective map. We set  $\xi(Z_0)=det(\eta_{Z_0}\circ \eta_{ 0}^{-1})$ and  $^\backprime \g z^-=\{Z\in\g z^-; \xi(Z)\neq 0\}$. Then  $^\backprime \g z^-$ is invariant under $H_s^+$.

Thus the map  $\ga$ from $ H\times^\backprime\g z^- $ to $\g q$ defined  by $\ga(h, Z) =h\cdot (A_0+Z)$ is a submersion. By Theorem \ref{rest}, 
 for any    $H$-invariant distribution  $\T$ on $\g q$, there exists a unique $H_s^+$-invariant distribution $\cR es_{\g z^-}\T$ defined on  $^\backprime \g z^-$  such that, for any $f\in\cD(H\times^{\backprime}\g z^-)$, one has  $<\T,\ga_{\star}(f)>=<\cR es_{\g z^-}\T,p_{\star}(f)>$. \me
 
Let   $\o_{\g z^-}$ be the  restriction of $\o$ to $\g z^-$. Then, one has:
 
 \begin{lemma}\label{radomega}(\cite{Se}) Lemma 4.4). Let  $\cR ad_{\g z^-}(\partial(\o))$ be the radial part of $\partial(\o)$ with respect to  $\ga$ (Theorem \ref{rad}). Then  $${\mathcal R}ad_{\g z^-}(\partial(\o))=\xi^{-1/2} \partial(\o_{\g z^-})\circ\xi^{1/2}-\mu$$ where  $\mu(Z)=\xi(Z)^{-1/2} \big(\partial(\o_{\g z^-})\xi^{1/2}\big)(Z)$ is an analytic function  on  $^\backprime \g z^-$. 
\end{lemma}

\section {Nilpotent and distinguished elements}\label{nilpotent}

Let $Z_0\in\g q$. Let $Z_0=A_0+X_0$ be its Jordan decomposition (\cite{Se} Lemma 1.1).  We construct the symmetric pair $(\g z_s,\g z_s^+)$ related to $A_0$ as in \ref{semisimple}. \me

We assume that $X_0 $ is different from zero. From (\cite{Se} Lemma 1.7), there exists  a  normal $sl_2$-triple  $(B_0, X_0, Y_0)$  of $(\g z_s,\g z_s^+)$ containing  $X_0$, i.e. satisfying $B_0\in \g z_s^+$ and $Y_0\in \g z_s^-$ such that $[B_0, X_0]=2 X_0$, $[B_0, Y_0]=-2 Y_0$ and $ [X_0, Y_0]=B_0$ . \me

We set  $\g z_0=\R B_0+\R X_0+\R Y_0$. The Cartan involution     $\t_0$ of $\g z_0$ defined by  $\t_0: (B_0,X_0,Y_0)\rightarrow (-B_0, -Y_0,-X_0)$ extends  to a Cartan involution of $\g z_s$, denoted by $\t$, which commutes with $\s$. (\cite{vd} Lemma 1). The bilinear form $(X,Y)\mapsto -B(\t(X),Y)$ defines a scalar product on $\g z_s$.\me

We can decompose $\g z_s$ in an orthogonal sum   $ \g z_s=\sum_{i}\g z_i$    of irreducible representations $\g z_i$ under the adjoint action of $\g z_0$. 
One can choose a suitable ordering  of the  $\g z_i$ such that  $  (\g z^-_s)_{Y_0}=\sum_{i=1}^r \g z_i\cap (\g z^-_s)_{Y_0}=\t((\g z^-_s)_{X_0})$
 with $\g z_1=\g z_0$  and $\textrm{dim } \g z_i\cap \g (\g z^-_s)_{Y_0}=1$. We set $n_i+1=\textrm{dim } \g z_i$. Hence, there exists an orthonormal basis  $(w_1,\ldots , w_r)$ of $(\g z^-_s)_{Y_0}$ such that $\dis w_1=\frac{Y_0}{\Vert Y_0\Vert}$  and  $[B_0, w_i]=-n_i w_i$  for $i\in\{1,\ldots, r\}$. In particular, one has $n_1=2$.  \me

We set
$$\delta_{\g q}(Z_0)=\delta_{\g z^-_s}(X_0)=\sum_{i=1}^r (n_i+2)-dim\;(\g z^-_s).$$

Let $\cN(\g z_s^-)$ be the set of nilpotent elements of $\g z_s^-$.

\begin{Definition}\label{defdistingue} (\cite{Se} Definitions 1.11 and 1.13)\begin{enumerate}
\item An element    $X_0$ of $\cN(\g z^-_s)$ is a $\g z^-_s$-distinguished nilpotent element    if  $ (\g z_s^-)_{X_0}$ contains no non-zero semisimple element.

\item An element  $Z_0$ of $\g q$ with Jordan decomposition $Z_0=A_0+X_0$ is called $\g q$-distinguished if  $X_0$ is a $\g z^-_s$-distinguished nilpotent element of  $\g z^-_s$.
\end{enumerate}\end{Definition}

\begin{Definition} \label{defSeki} The symmetric pair $(\g g,\g h)$ is nice if for any  $\g q$-distinguished element $Z$, one has $\delta_{\g q}(Z)>0$. 
\end{Definition}

Let  $\o_s$ be the restriction of $\o$ to $\g z^-_s$. Though $\o_s$ is not the Casimir polynomial on  $\g z_s^-$, one has the following result:

\begin{lemma} (\cite{vd} Lemma 4) The following assertions are equivalent:
\begin{enumerate}
\item $X_0$ is a $\g z^-_s$-distinguished nilpotent element.
\item $\o_s(X)=0$ for all $X\in (\g z_s^-)_{X_0}$.
\item  $\o_s(X)=0$ for all $X\in (\g z_s^-)_{Y_0}$.
\item $n_i>0$.
 \item $ (\g z_s^-)_{X_0}\cap  (\g z_s^-)_{Y_0}=\{0\}$.
 \end{enumerate}
 \end{lemma}
  
Thus, if $X_0$ is a $\g z^-_s$-distinguished nilpotent element  then  one has $\o(X_0+X)=2B(X_0,X)=2\Vert Y_0\Vert x_1$  for all  $X\in(\g z_s^-)_{Y_0}$, where $x_1$ is the first coordinate of $X$ in the basis $(w_1,\ldots , w_r)$ of $(\g z^-_s)_{Y_0}$.\\
  
  \no For  any $X_0\in\cN(\g z^-_s)$, 
one has
  $\g z_s^-=(\g z_s^-)_{Y_0}\oplus [\g z^+_s, X_0]$ and $\g z_s^+=(\g z_s^+)_{X_0}\oplus [\g z^-_s, Y_0].$
From now on, we set 
$$U=(\g z_s^-)_{Y_0}.$$
For $X\in U$, we consider the map $\psi_X$ from $[\g z^-_s, Y_0]\times U$ to $ \g z_s^-$ defined by $\psi_X(v,z)= z+[v,X_0+X]$.
The map $\psi_0$ is bijective. \\
We set  $\kappa(X)=det(\psi_X\circ \psi_0^{-1})$ and  $^\backprime  U=\{X\in U; \kappa(X)\neq 0\}.$ 
Hence, the map $\pi$ from $H_s^+\times ^\backprime  U$ to $ \g z_s^-$ defined by $\pi (h, X)=h\cdot(X_0+X) $ is a submersion.\me

We precise now some properties of $\pi$ related to  $\cN(\g z_s^-)$. \me

By (\cite{vd2} Theorem 23]), we can write  $\cN(\g z_s^-)=\cO_1\cup\ldots\cO_\nu$ where the $\cO_j$ are disjoints $H_s^+$-orbits with $\cO_\nu=\{0\}$  and each  $\cO_j$  is open in the closed set $\cN_j=\cO_j\cup\ldots\cO_\nu$. 
One assumes that $\cO_j=H_s^+\cdot X_0$.

\begin{lemma} \label{choixU0}(\cite{vd} Lemma 17 and 18). There exists a neighborhood  $U_0$ of $0$ in $ U$ such that
\begin{enumerate} \item $\pi$ is a submersion on $H_s^+\times U_0$,
\item $\O_0=\pi(H^+_s\times U_0)$ is an open neighborhood of $X_0$ in $\g z_s^-$ and $\O_0\cap\cN_j=\cO_j$,
\item $\cO_j\cap (X_0+U_0)=\{X_0\}$
\item Let $\T$ be an  $H_s^+$-invariant distribution  on $\O_0$. Let  $\cR es_U \T$ be its restriction to $U$ with respect to  $\pi$.

If $\textrm{ supp }(\T)\subset \cN_j$ then $\textrm{ supp }({\mathcal R}\mathit{es}_U \T)\subset\{0\}$.
\end{enumerate}
\end{lemma}

We denote by $\o_{\g c^-}$ and  $\o_s$  the restrictions of $\o$   to $\g c^-$ and $\g z_s^-$  respectively. One has   $\o_{\g z^-}= \o_{\g c^-}+\o_s$.   We precise now the radial part  ${\cal R}ad_U(\partial(\o_s))$ of   $\partial(\o_s)$  with respect to $\pi$. We denote by   ${\cal R}ad_{U,X}(\partial(\o_s))$ its local expression at $X\in U_0$.

\begin{lemma}\label{nondistingue}(\cite{vd} Lemma 13) 
 The homogeneous part of degree 2 of    ${\cal R}ad_{U,0}(\partial(\o_s))$ is zero if and only if  $X_0$ is   $\g z_s^-$-distinguished.
 \end{lemma}
 \begin{Theorem}\label{distingue}(\cite{vd}  Theorem  14)
Let   $X_0$ be a $\g z_s^-$-distinguished nilpotent element and $c_0=\Vert X_0\Vert$. Then, there exist analytic functions $a_{i,j}$ ($2\leq i,j\leq r$) and $a_i$ ($2\leq i\leq r$) on $U_0$ satisfying  $a_{i,j}(0)=0$ such that, for any  $H_s^+$-invariant distribution $T$ on  $\O_0$, one has

$$  {\mathcal R}es_U(\partial(\o_s)T)=  {\cal R}ad_U((\partial(\o_s)){\mathcal R}es_U(T)$$
$$=\frac{1}{c_0}\Big(2x_1\frac{\partial^2}{\partial x_1^2}
+(dim\; \g z_s^-) \frac{\partial}{\partial x_1}+\sum_{i=2}^r(n_i+2)x_i\frac{\partial^2}{\partial x_1\partial x_i}$$
$$
+\sum_{2\leq i\leq j\leq r} a_{i,j}(X) \frac{\partial^2}{\partial x_j\partial x_i}+\sum_{i=2}^r a_i(X)\frac{\partial}{\partial x_i}\Big){\mathcal R}es_U(T)$$
where   $x_1,\ldots,x_r$ are the coordinates of $X$ in the basis  $(w_1,\ldots , w_r)$.\end{Theorem}
 
 \section{ The main Theorem}
Our goal is to prove the following Theorem:
\begin{Theorem}\label{resultatprincipal} Let $(\g g, \g h)$ be a nice reductive symmetric pair. Let $\cV$ an $H$- invariant open subset of $\g q$. Let $\T$ be an $H$-invariant distribution on  $\cV$ such that
\begin{enumerate}\item There exists  $P\in\C[X]$ such that $P\big(\partial(\o)\big) \T=0$
\item There exists $F\in L^1_{loc}(\cV)^H$ such that  $\T=F$ on  $\cV\cap\g q^{reg}$.
\end{enumerate}
Then $\T=F$ as distribution on  $\cV$. \end{Theorem}

We will use the method developed by M.  Atiyah  in \cite{A}. First we recall some facts about  distributions on $\R^r   \times \R^m$.  Let $\N$ be the set of non-negative integers. For $\a=(\a_1, \ldots,\a_r)\in\N^r$, we set $\vert \a\vert=\a_1+\ldots +\a_r$ and
$$x^\a=x_1^{\a_1}\ldots x_r^{\a_r}, \quad \partial_x^\a=\frac{\partial^{\vert\a\vert}}{\partial x_1^{\a_1}\ldots \partial x_r^{\a_r}}.$$
 
For  $\varphi\in \cD(\R^r \times\R^{m})$ and $\e>0$, we set $\varphi_\e(x, y)=\varphi(\frac{x}{\e}, y)$ for $(x, y)\in \R^r\times\R^m$. For  $T\in\cD'(\R^r \times \R^m)$ we denote by $T_\e$ the distribution defined by
$\displaystyle <T_\e,\varphi >=<T,\varphi_\e >$. 

\begin{Definition} Let $V=\{0\}\times \R^m\subset \R^r \times\R^m$ and $T\in\cD'(\R^r  \times\R^m)$.
\begin{enumerate}\item The distribution $T$ is regular  along $V$ if $\displaystyle \lim_{\epsilon\rightarrow 0} T_\epsilon=0$. 

\item The distribution $T$ has a   degree  of singularity    along $V$ smaller than $k$  if for all $\a\in\N^r$ with $\vert\a\vert=k$, the distribution $x^\a T$ is regular.

We denote by $d^{\circ}_s T$ the degree of singularity  of $T$  along $V$ and we omit in what follows to precise "along $V$". Regularity corresponds to a  degree of singularity   equal to $0$.
\item  The degree of singularity  of $T$ is equal to $k$ if $d^{\circ}_s T\leq k$ and $d^{\circ}_s T\nleq k-1$.

\end{enumerate}
\end{Definition}

\begin{lemma}\label{degres}  
\begin{enumerate}

\item If $F\in L^1_{loc}( \R^{r+m})$ then $d^\circ_s F=0$.
\item If $d^\circ_s T= k\geq 1$ then $d^\circ_s (x_iT)= k-1$ for $i\in\{1,\ldots r\}$.
\item  If $d^\circ_s T\leq k$ then $\displaystyle\frac{\partial}{\partial x_i} T\leq k+1$ for $i\in\{1,\ldots r\}$.
\item Let  $\delta_0$  be the  Dirac measure at $0\in\R^r $ and  $\delta_0^{(\a)}= \partial_x^{\a}\delta_0$. If $S\in\cD'(\R^m)$ then the degree of singularity of  $\delta_0^{(\a)}\otimes S$ is equal to $\vert\a\vert+1$.
\end{enumerate}
\end{lemma}
\no{\it Proof.} {\it 1.} Let $F\in L^1_{loc}( \R^{r+m})$  and $\phi\in\cD( \R^{r+m})$ with $\mbox{supp}(\phi)\subset K_1\times K_2$ where $K_1$ (resp., $K_2$) is a compact subset of $\R^r$ (resp., $\R^m$). One has 
$$\vert\int_{\R^r\times \R^m} F(x,y) \phi(\frac{x}{\e},y) dx dy\vert\leq \sup_{(x,y)\in \R^{r+m}}\vert \phi(x,y)\vert\int_{(\e K_1)\times K_2}\vert F(x,y)\vert dx dy$$
and  the first assertion follows.

\no {\it 2.} is clear. 

\no {\it 3.} Let $\a\in\N^n$ such that $\vert \a\vert=k+1$. If $\a_j\geq 1$ for some $j\in\{1,\ldots,r\}$, we set $\bar{\a}^j=(\a_1,\ldots,\a_{j-1},\a_j-1,\a_{j+1},\ldots, \a_r)$. Let $\varphi\in\cD(\R^{r+m})$. 

If $\a_i\geq 1$, one has $$<x^\a\frac{\partial}{\partial x_i} T,\varphi_\e>=- <T,\a_ix^{\bar{\a}^i}\varphi_\e +\frac{x^\a}{\e}(\frac{\partial}{\partial x_i}\varphi)_\e>$$
$$=-\a_i<x^{\bar{\a}^i}T, \varphi_\e>-<x^{\bar{\a}^i}T, (x_i\frac{\partial}{\partial x_i}\varphi)_\e>$$
thus  $(x^\a T)_\e$ converges to $0$ 
since $d^\circ_sT\leq k$.

If $\a_i=0$, we choose $j$ such that $\a_j\geq 1$. One has $\displaystyle <x^\a\frac{\partial}{\partial x_i} T,\varphi_\e>=-<x^{\bar{\a}^j}T, (x_j\frac{\partial}{\partial x_i}\varphi)_\e>$ which tends to $0$ as before.

\no{\it 4.} We recall that for $i\in\{1,\ldots,r\}$, one has

$$x_i^l\delta_0^{(\a)}=\left\{\begin{array}{cc} (-1)^l\frac{(\a_i)!}{(\a_i-l)!}\delta_0^{(\a_1,\ldots, \a_i-l,\ldots \a_n)}& \mbox{ if }\a_i\geq l\\

0& \mbox{ if }\a_i<l.\end{array}\right.$$
Hence, one has $x^\a\delta_0^{(\a)}=(-1)^{\vert\a\vert}\a!\delta_0$ and for all $\b\in\N^r$ with $\vert \beta\vert=\vert\a\vert+1$, one has $x^\b\delta_0^{(\a)}=0$. The assertion follows.\qed

\begin{Definition} Let $\G=x^\beta\partial^\a_x D$  where $D$ is a differential operator on $\R^m$. Then $\G$ increases the  degree of singularity   at most $\vert\a\vert-\vert\b\vert$. The integer $\vert\a\vert-\vert\beta\vert$ is called the total degree of   $\G$ in $x$.

We can define the homogeneous part of highest total degree (in $x$) of an analytic differential operator developing its coefficients in Taylor series.\end{Definition}

\no{\bf Proof of the Theorem.} Let  $\T\in\cD'(\cV)^{H}$ and $F\in L_{loc}^1(\cV)^H$ such that
 $P(\partial(\o)) \T=0$ for a unitary polynomial $P\in\C[X]$  and    $\T=F$ on $\cV^{reg}=\cV\cap\g q^{reg}$. 
We write    $\T=F+S$ where $S$ is an  $H$-invariant distribution with support contained in  $\cV-\cV^{reg}$. We want to prove that  $S=0$, which is  equivalent to $\textrm{ supp }(S)=\emptyset$. \me

Assuming  $S $ is non-zero, we are led to a contradiction. We will study  $S$ near an element  $Z_0\in \textrm{ supp } (S)$ chosen as follows:

For   $Z_0\in \textrm{ supp }(S)$ with  Jordan decomposition $Z_0=A_0+X_0$, we construct the symmetric subpair $(\g z_s,\g z_s^+)$  related to $A_0$ and we set   $\g q_{A_0}= \g z^-=\g c^-\oplus \g z_s^-$ as in section  \ref{semisimple}. Let $\cS_k$ be the set of $Z_0$ in the support of $S$ such that $\textrm{rank}(\g z_s^-)=k $. Since$\textrm{ supp } (S)\subset \cV-\cV^{reg}$, if $Z_0=A_0+X_0$ belongs to  $\textrm{ supp }(S)$  then $A_0$ is not $\g q$-regular. One deduces that $S_0=\emptyset$. Let  $k_0>0$ such that  $S_0=S_1=\ldots =\cS_{k_0-1}=\emptyset$ and  $\cS_{k_0}\neq \emptyset$. 

For   $Z_0=A_0+X_0$ in $ \cS_{k_0}$, we denote by   $\cN(\g z_s^-)=\cO_1\cup\ldots\cO_\nu$ the set of nilpotent elements in $\g z_s^-$ as in section \ref{nilpotent}. Since $\textrm{ supp } (S)\cap (A_0+\cN(\g z_s^-))\neq\emptyset$, one can choose
 $j_0\in \{1,\ldots,\nu \}$  such that  $\textrm{ supp }(S)\cap (A_0+\cO_{i})=\emptyset$ for $i\in\{1,\ldots j_0-1\}$ and $\textrm{ supp }(S)\cap (A_0+\cO_{j_0})\neq\emptyset$.  \me

From now on, we fix  $Z_0=A_0+X_0 $ in 
$ \cS_{k_0}$ such that $X_0\in\cO_{j_0}$. \me

For $\e>0$, we denote by  $\cW_\e$   the set of $x$ in $\g z_s^-$ such that, for any eigenvalue $\l$ of $\textrm{ad}_{\g g} x$, one has $\vert\l\vert<\e$.
The choice of $k_0$ implies that there exists $\e>0$ such that $ \mbox{ supp}(S)\cap (Z_0+\cW_\e)\subset  \mbox{ supp}(S)\cap(Z_0+\g c^-+\cN(\g z_s^-))$. Hence, we can choose an open neighborhood  $\cW_c$ of $0$ in $\g c^-$ and an open neighborhood $\cW_s$ of $X_0$ in $\g z_s^-$ such that

\begin{equation}\label{ouvert} \mbox{ supp}(S)\cap (A_0+\cW_c+\cW_s)\subset \textrm{ supp }(S)\cap(A_0+\cW_c+\cN(\g z_s^-)).\end{equation}

       \no{\bf First case.} $A_0\notin \g c_{\g q}$ and $X_0\neq 0$. 

We keep the notation of section  \ref{nilpotent}. We fix a normal  $sl_2$-triple $(B_0,Y_0,X_0)$ in $(\g z_s,\g z_s^+)$. We choose   an open neighborhood $U_0$ of $0$ in $U$, the centralizer of  $Y_0$ in $\g z_s^-$,  as in Lemma   \ref{choixU0}. We keep the notation of this lemma. We recall that   the map $\ga$ from $ H\times^\backprime\g z^- $ to $ \g q$ defined by $\ga(h, Z) = h\cdot (A_0+Z) $ is a submersion. Reducing  $U_0$, $\cW_c$ and $\cW_s$ if necessary, we may assume that $ \cW_c+\O_0\subset \cW_c+\cW_s\subset ^\backprime \g z^-$ and  that $V_0=\ga(H\times (\cW_c+\O_0))$  is an open neighborhood of  $Z_0$ contained in $\cV$. \me

If $T$ is an  $H$-invariant  distribution on $\cV$, we denote by $T_0$ its restriction to $V_0$. By theorem \ref{rest}, one can consider its restriction  $T_1={\mathcal R}es_{\g z^-} T_0$ to $\cW_c+\O_0$ with respect to $\ga$. One has $A_0+\textrm{ supp }(T_1)\subset \textrm{ supp } (T)\cap(A_0+ \cW_c+\O_0)$.

We set $T_2=\xi^{1/2} T_1$ where $\xi^{1/2}$ is the analytic function on $\cW_c+\O_0$ defined in    section \ref{semisimple}.

Now, we consider  the  submersion $\pi_0$ from $ H_s^+\times U_0\times\cW_c$ to  $\g z^-$ defined by $\pi_0 (h, X,C)= h\cdot(X_0+X)+C $. One denotes by $T_3 $ the restriction on $ U_0\times \cW_c$ of $T_2$ with respect to $\pi_0$ . We have $X_0+\textrm{supp}(T_3)\subset \textrm{ supp} (T_2)\cap(X_0+U_0)$.\me

Since $F$ is a locally  integrable function, the distribution  $F_3$ is the locally  integrable function  on  $U_0\times \cW_c$ defined by  $F_3(X,C)=\xi^{1/2}(C+X) F(C+X)$.\me

By assumption, the distribution $S_3$ is non-zero.  By $(\ref{ouvert})$ and Lemma \ref{choixU0} ({\it 2.}), one has $\textrm{ supp }(S_2)=\textrm{ supp }(S_1)\subset \cW_c+\O_0\cap\cN_{j_0}=\cW_c+\cO_{j_0}$. We deduce from Lemma \ref{choixU0} ({\it 3.}) that  $\textrm{ supp }(S_3)\subset \{0\}\times \cW_c$. By (\cite{Sc},  Lemma 3), there exists a family $(S_\a)_\a$ of  $\cD'(\cW_c)$ such that  
$\dis S_3=\sum_{\a\in\N^r; \vert\a\vert  \leq l}  \delta^{(\a)}_{0}\otimes S_\a$ where $\delta_0$ is the Dirac measure at  $0$ of $U_0$ and for $\a\in\N^r$ , the $S_\a$  with $\vert \a\vert=l$ are not all zero. \me

By assumption, the distribution $\T$ satisfies $P\big(\partial(\o) \big)\T=0$. By Lemma \ref{radomega}, one has

$$P\Big((\partial(\o_{s})+\partial(\o_{\g c}) )-\mu(Z)\Big)\T_2=0 \textrm{ on } \cW_c+\O_0.$$

\no Using the restriction with respect to $\pi_0$, one obtains
$$P\Big(\cR ad_U(\partial(\o_s))+\partial(\o_{\g c})-\tilde{\mu}\Big)\T_3=0 \textrm{ on   }U_0 \times   \cW_c$$
where $\tilde{\mu}(X,C)=\mu(C+X)$ for $X\in U_0$ and $C\in \cW_c$.\me

 Let $D_0$ be the homogeneous part of highest total degree $d$ of $\cR ad_{U }(\partial(\o_s))$.
 We set
$$P\Big(\cR ad_U(\partial(\o_s))+\partial(\o_{\g c})-\tilde{\mu}\Big)=D_0^N+D_1$$ where    $N$ is the degree of $P$ and $D_1$ is a differential operator with total degree in  $X $ strictly smaller than $Nd$.
Since $\T_3=F_3+S_3$ with $\dis S_3=\sum_{\\a\in\N^r; \a_1  \leq l}  \delta^{(\a)}_{0}\otimes S_\a$, we obtain the following relation on $U_0\times \cW_ c$:
\begin{equation}\label{etude}(D_0^N+D_1) S_3=(D_0^N+D_1) ( \sum_{\a\in\N^r; \vert \a\vert\leq l}  \delta^{(\a)}_{0}\otimes S_\a)= -(D_0^N+D_1) F_3 \end{equation}

We study now  the degree of singularity along $\{0\}\times\cW_c$ of the two members of  $(\ref{etude})$.  \me

If $X_0$ is not a $\g z_s^-$-distinguished nilpotent element then by  Lemma \ref{nondistingue}, the homogeneous part of degree $2$ of $\cR ad_{U,0}(\partial(\o_s)$ does not vanish and is a differential operator with constant coefficients of degree $2$. Hence  the total degree   of $D_0$ is  equal to  $d=2$.   Since $F_3 $ is a locally integrable function, it follows from Lemma \ref{degres} that one has $d^\circ_s F_3=0$ and $d^\circ_s ((D_0^N+D_1) F_3)\leq 2N$. By the same Lemma, one has $d^\circ_s ((D_0^N+D_1) S_3)= l+1+2N$. Hence,  we have  a contradiction.\me

Assume that   $X_0$ is  a $\g z_s^-$-distinguished nilpotent element.  Lemma \ref{distingue} gives   $\dis  c_0 D_0= 2 x_1\frac{\partial^2}{\partial x_1^2}
+(dim\; \g z_s^-) \frac{\partial}{\partial x_1}+\sum_{i=2}^r(n_i+2)x_i\frac{\partial^2}{\partial x_1\partial x_i}+\sum_{2\leq i\leq j\leq r} a_{i,j}(X) \frac{\partial^2}{\partial x_j\partial x_i}+\sum_{i=2}^r a_i(X)\frac{\partial}{\partial x_i})$  where    $c_0=\Vert X_0\Vert$ . Since $a_{i,j}(0)=0$, the   total degree of $D_0$  is equal to $1$. \me

For $\a=(\a_1,\ldots,\a_r)\in\N^r$, we set $\tilde{\a}^i=(\a_1,\ldots,\a_{i-1},\a_i+1,\a_{i+1}\ldots \a_r)$ and $\bar{\a}^i=(\a_1,\ldots,\a_{i-1},\a_i-1,\a_{i+1}\ldots \a_r)$. The relation $x_i\delta_0^{(\a)}=-\a_i\delta_0^{(\bar{\a}^i)}$ and the above expression of $D_0$ give
$$c_0 D_0\cdot \delta_0^{(\a)}\otimes S_\a=\l_\a \delta^{(\tilde{\a}^1)}\otimes S_\a+\sum_{2\leq i\leq j\leq r} a_{i,j}(X) \delta^{(\tilde{\a}^{i,j})}\otimes S_\a+\sum_{i=2}^r a_i(X)\delta^{(\tilde{\a}^i)}\otimes S_\a$$ where   $$  \l_\a= -2(\a_1+2)+dim\; \g z_s^- -\sum_{i=2}^r (n_i+2)(\a_i+1).$$

Since $n_1$ is equal to $2$ and  $(\g g, \g h)$ is a nice pair, we obtain

$$
 \l_\a=-\delta_{\g q}(Z_0)-\big[ 2\a_1+\sum_{i=2}^r (n_i+2)\a_i\big ]<0\textrm{  for all  } \a\in\N^r.$$

 Consider  $\a_0=(\a_1,\ldots,\a_r)\in\N^r$ such that $\vert\a_0\vert=l$, $S_{\a_0}\neq 0$ and $\a_1$ is maximal for these properties. One deduces that the coefficient of $\delta^{(\widetilde{\a_0}^1)}\otimes S_{\a_0}$ in $D_0\cdot(\sum_{\a\in\N^r;\vert\a\vert=l} \delta_0^{(\a)}\otimes S_\a)$ is non-zero. Thus, the degree of singularity of $(D_0^N+D_1)S_3$ is equal to $1+l+N$. Since $F_3$ is locally integrable and the  total degree of $D_0$ is equal to $1$, we have $d^\circ_s(D_0^N+D_1)F_3\leq N$. This gives a contradiction in (\ref{etude})\me

\no{\bf Second  case.}   $A_0\in\g c_{\g q}$ and $X_0\neq 0$. 

The symmetric pair $(\g z_s,\g z_s^+)$ is equal to $(\g g_s, \g h_s)$.  We just consider the submersion 
$\pi_0$ from $ H\times U_0\times\cW_c $ to $ \g q$ defined by $\pi_0 (h, X,C)=h\cdot(X_0+X)+A_0+C$
where $U_0$ is defined as in Lemma \ref{choixU0} for the symmetric pair $(\g g_s, \g h_s)$.

For $T\in\cD'(\g q)^H$, we denote by  $T_1$   the restriction of $T$ to $U_0\times\cW_c$ with respect to $\pi_0$. As in the first case, we have  $\T_1=F_1+S_1$ where $F_1$ is a locally integrable function on  $U_0 \times   \cW_c$ and  $S_1$ is a non-zero distribution such that $\textrm{ supp }(S_1)\subset \{0 \}\times   \cW_c$. Moreover the  distribution  $\T_1$ satisfies the relation

$$P\Big(\cR ad_U(\partial(\o_s))+\partial(\o_{\g c})\Big)\T_1=0 \textrm{ on   }U_0 \times   \cW_c.$$

\no The same arguments as in the first case lead to  the contradiction $S_1=0$.\me

\no{\bf Third case.}    $X_0= 0$. 

The open sets $\cW_c$ and $\cW_s$ satisfy $\textrm{ supp }(S)\cap(A_0+\cW_c+\cW_s)\subset \textrm{ supp }(S)\cap (A_0+\cW_c+\cN(\g z_s^-))$. By the choice of $j_0$, we deduce that  $\textrm{ supp }(S)\cap(A_0+\cW_c+\cW_s)\subset \textrm{ supp }(S)\cap (A_0+\cW_c)$. \me

If $A_0\in\g c_{\g q}$, then $V_0=A_0+\cW_c+\cW_s$ is an open neighborhood of $A_0$ in $\g q$. We identify $\g q$ with $\g q_s\times \g c_{\g q}$. Thus, the restriction $S_0$ of $S$ to $V_0$ is different from zero  and satisfies    $\textrm{supp} (S_0)\subset\{0\}\times (A_0+\cW_c)$. On the other hand, one has  $P(\partial(\omega) )S_0=-P(\partial(\omega) )F_{|V_0}$. Since $\partial(\o)$ is a second order operator with constant coefficients,   we obtain a contradiction as above.\me

If $A_0\notin\g c_{\g q}$, we may assume that  $\cW_c+\cW_s\subset ^\backprime\g z^-$. We denote by  $T_1$ the restriction of an $H$-invariant distribution  $T$ to  $\cW_c+\cW_s$ with respect to the submersion $\ga$ from $H\times ^\backprime\g z^-  $ to $  \g q $ and  we consider $T_2=\xi^{1/2}T_1$ as   distribution on $\cW_s\times\cW_c$. Thus, we have    $S_2\neq 0$ and $\textrm{ supp } (S_2)=\{0\}\times  \cW_c$. Moreover, the distribution  $\T_2=F_2+S_2$ satisfies 
$P\Big((\partial(\o_{s})+\partial(\o_{\g c}) )-\mu(Z)\Big)\T_2=0$ on $ \cW_s\times \cW_s$ by Lemma \ref{radomega}. This is equivalent to
$$P\Big((\partial(\o_{s})+\partial(\o_{\g c}) )-\mu(Z)\Big)S_2=-P\Big((\partial(\o_{s})+\partial(\o_{\g c}) )-\mu(Z)\Big)F_2.$$
Since  $\partial(\o_s)$ is a second order operator with constant coefficients,   we obtain a contradiction as above.

This achieves the proof of the Theorem.\qed
\section{ Application to $\big(\mathfrak{gl}(4,\R), \mathfrak{gl}(2,\R)\times \mathfrak{gl}(2,\R)\big)$}
On  $G=GL(4,\R)$ and its Lie algebra  $\g g=\mathfrak{gl}(4,\Bbb R)$, we consider the involution $\sigma$ defined by $\displaystyle \sigma(X)=\left
(\begin{array}{cc}
I_{2}&0\\0&-I_{2}
\end{array}\right )X
\left (
\begin{array}{cc}I_{2}&0\\0&-I_{2}
\end{array}\right )$
where $I_{2}$ is the  $2\times 2$ identity matrix. 
 We have 
 $\g g=\g h\oplus\g q$ with $$\g h=\left\{ \left(
\begin{array}{cc} A& 0\\ 0&B
\end{array}\right ); A, B\in  {\mathfrak gl}(2,\R)\right\} \textrm{  and } \g q=\left\{ \left (
\begin{array}{cc}0&Y\\Z&0
\end{array}\right );  Y,Z\in   {\mathfrak gl}((2,\R)\right\}.$$ By (\cite{Se} Theorem 6.3), the symmetric pair $\big( {\mathfrak gl}((4,\R), {\mathfrak gl}((2,\R)\times {\mathfrak gl}((2,\R)\big)$ is a nice pair.

We first recall some results of \cite{HJ}.
Let $\kappa(X,X')=\dfrac{1}{2} tr(XX')$. The restriction of $\kappa$ to the derived algebra of $\g g$ is a multiple of the Killing form. Let $S(\g q_\C)^{H_\C}$ be subalgebra  of  $S(\g q_\C)$ of all elements invariant under $H_\C$. We identify $S(\g q_\C)^{H_\C}$ with the algebra of $H_\C$-invariant differential operators on $\g q_\C$ with constant coefficients. Using  $\kappa$, we identify  $S(\g q_\C)^{H_\C}$ with the algebra $\C[\g q_\C]^{H_\C}$ of $H_\C$-invariant polynomials  on $\g q_\C$. A basis of  $\C[\g q_\C]^{H_\C}$ is given by 
   $Q(X)=\dfrac{1}{2} tr (X^2)$ and $S(X)=det(X)$. The Casimir polynomial is just a multiple of $Q$.
   
By (\cite{HJ} Lemma 1.3.1), the $H$-orbit of a semisimple element $X=\left(\begin{array}{cc}0&Y\\Z&0
\end{array}\right)$ of $\g q$ is characterized by $(Q(X),S(X))$ or by the set $\{ \nu_1(X), \nu_2(X)\}$ of eigenvalues of $YZ$,  where the functions $\nu_1$ and $\nu_2$ are defined as follows: let $Y$ be the Heaviside function. Let $S_0 =Q^2 -4S$ and $\d=\iota^{Y(-S_0)}\sqrt{\vert S_0\vert}$. We set 
$$\nu_1=(Q+\d)/2\quad \textrm{ and }\quad \nu_2=(Q-\d)/2.$$ Regular elements of $\g q$ are semisimple elements with 2 by 2 distinct eigenvalues or equivalently, semisimple elements $X$ of $\g q$ such that $ \nu_1(X) \nu_2(X)( \nu_1(X)- \nu_2(X))\neq 0$ (\cite{HJ} Remarque 1.3.1). \me

Let $\chi$ be the character of $\C[\g q_\C]^{H_\C}$ defined by $\chi(Q)=\l_1+\l_2$ and $\chi(S)=\l_1\l_2$ where  $\l_1$ and $\l_2$ are two complex numbers satisfying $\l_1\l_2 (\l_1-\l_2)\neq 0$. 

For an open $H$-invariant subset $\cV$ in $\g q$, we denote by $\cD' (\cV)^H_{\chi}$ the set   of $H$-invariant distributions $T$ with support in  $\cV$ such that $\partial(P) T=\chi(P) T$ for all $P\in\C[\g q_\C]^{H_\C}$.  
Let $\cN$ be the set of nilpotent elements of $\g q$ and $\cU=\g q-\cN$   its complement. In \cite{HJ}, we describe a basis of the subspace of $ \cD' (\cU)^H_{\chi}$ consisting of locally integrable functions. More precisely, we obtain the following result.\me

We consider the Bessel operator $ L_c=4\left(z\frac{\partial^2}{\partial z^2}+\frac{\partial}{\partial z}\right)$ on $\C$ and its analogous $L=4\left(t\frac{d^2}{dt^2}+\frac{d}{dt}\right) $ on $\R$. Let  $\mathcal{S}ol(L_c,\l)$ (resp., $\mathcal{S}ol(L ,\l)$) be the set of holomorphic (resp., real analytic ) functions $f$ on $\C-\R_-$ (resp., $\R^*$)  such that $L_c f=\l f$ (resp., $L  f=\l f$). 
For $\l\in\C^*$, we set $$\Phi_\l(z)=\sum_{n\geq0}\frac{(\l z)^n}{4^n(n!)^2}\quad \textrm{ and }\quad w_\l(z)=\sum_{n\geq0}\frac{a(n)(\l z)^n}{4^n(n!)^2},$$  where $ a(x)= -2 \frac{\Gamma'(x+1)}{\Gamma(x+1)}.$
Then $(\Phi_\l, W_\l=w_\l+\log(\cdot)\Phi_\l)$ form a basis of $\mathcal{S}ol(L_c,\l)$ , where $\log$  is the principal determination of the 
  logarithm function on $\C-\R_-$ and $(\Phi_\l, W^r_\l=w_\l+\log\vert \cdot\vert \Phi_\l)$ form  a basis of $\mathcal{S}ol(L,\l)$.\me

For two functions $f$ and $g$ defined over $\C$, we set
 $$S^+(f,g)(X)=f( \nu_1(X))g( \nu_2(X))+f( \nu_2(X))g( \nu_1(X))$$
and  $$[f,g](X)=f( \nu_1(X))g( \nu_2(X))-f( \nu_2(X))g( \nu_1(X)).$$
We define the following functions on $\g q^{reg}$: 
\begin{enumerate}
    \item  $$F_{ana}=\frac{[\Phi_{\l_1},\Phi_{\l_2}]}{ \nu_1- \nu_2}$$

    \item  $$F_{sing}=\frac{[\Phi_{\l_1},w_{\l_2}]+[w_{\l_1},\Phi_{\l_2}]+\log|\nu_1\nu_2|[\Phi_{\l_1},\Phi_{\l_2}]}{\nu_1- \nu_2}$$
    \item   For 
$(A,B)\in\{(\Phi_{\l_1},\Phi_{\l_2}),\;(\Phi_{\l_1},W^r_{\l_2}),\;(W^r_{\l_1},\Phi_{\l_2}),\;(W^r_{\l_1},W^r_{\l_2})\}$, we set   $$F^+_{A,B}=Y(S_0)\dfrac{S^{+}(A,B)}{|\nu_1- \nu_2|}$$ where $S_0=Q^2-4S\in\C[\g q_\C]^{H_\C}$ and $Y$ is the Heveaside function.
\end{enumerate}
\begin{Theorem}\label{resultatHJ} (\cite{HJ} Theorem 5.2.2 and Corollary 5.3.1).
\begin{enumerate}\item The functions  $F_{ana}$ and $F_{sing}$ 
are locally integrable on  ${\g q}$.
 \item  
For 
$(A,B)\in\{(\Phi_{\l_1},\Phi_{\l_2}),\;(\Phi_{\l_1},W^r_{\l_2}),\;(W^r_{\l_1},\Phi_{\l_2}),\;(W^r_{\l_1},W^r_{\l_2})\}$, the functions  $F^+_{A,B}$, are locally integrable on  ${\mathcal U}$.
\item The family  $F_{ana}$, $F_{sing}$ and  $F^+_{A,B}$,
with $(A,B)$ as above 
  form a basis $\cB$ of the subspace of  $ \cD' (\cU)^H_{\chi}$ consisting of   distributions  given by a locally integrable function.
\end{enumerate}
\end{Theorem}

\begin{cor}\label{surU} Any invariant distribution of $ \cD' (\cU)^H_{\chi}$ is given by a locally integrable function on $\cU$. In particular, the family $\cB$ defined in the previous Theorem is  a basis of  $ \cD' (\cU)^H_{\chi}$.

\end{cor}
 
 \no{\it Proof.} Let  $T\in \cD' (\cU)^H_{\chi}$. We denote by  $F$ its restriction   to $\cU^{reg}$. By (\cite{Se} Theorem 5.3  (i)), $F$ is an analytic function on $\cU^{reg}$ satisfying $(*)\quad \partial(P) F= \chi(P) F$ on $\cU^{reg}$  for all $P\in\C[\g q _\C]^{H_\C}.$
 
 In (\cite{HJ} section 4.), we describe the analytic solutions of $(*)$ in terms of $\Phi_\l$, $W_\l$ and $W_\l^r$ for $\l=\l_1$ and $\l_2$. By the asymptotic behaviour of orbital integrals near non-zero semisimple elements (\cite{HJ} Theorems 3.3.1 and 3.4.1), and the Weyl integration formula (\cite{HJ} Lemma 3.1.2), one deduces that $F\in L^1_{loc}(\cU)^H$. 
 Theorem \ref{resultatprincipal} gives the result. \qed
 \me
 
\begin{cor} Any invariant distribution of $ \cD' (\g q)^H_{\chi}$ is given by a locally integrable function on $\g q$.
\end{cor}

 \no{\it Proof.}   Let  $T\in \cD' (\g q)^H_{\chi}$. By Corollary \ref{surU}, the restriction of $T$ to $\cU$ is a linear combination of elements of $\cB$. By Theorem \ref{resultatprincipal} and Theorem \ref{resultatHJ}, it is enough to prove that the functions $F^+_{A,B}$, with  $(A,B)\in \{(\Phi_{\l_1},\Phi_{\l_2}),\;(\Phi_{\l_1},W^r_{\l_2}),\;(W^r_{\l_1},\Phi_{\l_2}),\;(W^r_{\l_1},W^r_{\l_2})\}$ are locally integrable on $\g q$ or equivalently, that  the integral $\displaystyle \int_{\g q}| F^+_{A,B} (X)  f(X) | dX$ is finite for all positive function  $f\in\cD(\g q)$. For this, we will use the Weyl integration formula (\cite{or} Proposition 1.8 and Theorem 1.27).\me
 
For $\e=(\e_1,\e_2)$ with $\e_j= \pm$, we define  $$\g a_{\e }=\left\{X_{ \e}(u_1,u_2)=\left(\begin{array}{c|c} 0 & \begin{array}{cc} u_1 & 0\\   0& u_2 \end{array}\\\hline\begin{array}{cc}  \e_1 u_1 & 0\\  0&  \e_2 u_2 \end{array}&0\end{array}\right); (u_1,u_2)\in\R^2\right\}.$$

and $$\g a_2=\left\{ \left(\begin{array}{c|c}
0&\begin{array}{cc} \tau&-\t\\ \t&\tau\end{array}\\ \hline
\begin{array}{cc} \tau&-\t\\ \t&\tau\end{array}&0
\end{array}
\right); (\t,\tau)\in\R^2\right\}$$
By (\cite{HJ}, Lemma 1.2.1), the subspaces $\g a_{++},\g a_{+-}, \g a_{--}$  and $\g a_2$ form a system  of representatives    of $H$-conjugaison classes of Cartan subspaces in $\g q$.  By (\cite{HJ} Remark 1.3.1), an element $X\in \g q$ satisfies $S_0(X)\geq 0$ if and only if $X$ is $H$-conjugate to an element of $\g a_\e $ for some $\e$. Furthermore, one has $\{\nu_1(X_\e(u_1,u_2)), \nu_2(X_\e(u_1,u_2))\}= \{\e_1u_1^2, \e_2 u_2^2\}$. \me

Let $f$ be a positive function in $\cD(\g q)$. We define the orbital integral of $f$ on $\g q^{reg}$ by
$$\cM(f)(X)=|\nu_1(X)-\nu_2(X)|\int_{H/Z_H(X)} f(h.X) dX$$
where $Z_H(X)$ is the centralizer of $X$ in $H$ and $dh$ is an invariant measure on $H/Z_H(X)$. \me

By (\cite{or} Theorem 1.23), the orbital integral $\cM(f)$ is a smooth function on $\g q^{reg}$ and there exists a compact subset $\O$ of $\g q$ such that  $\cM(f)(X)= 0$ for all regular element $X$ in the complement of $\O$. \me

Since $F^+_{A,B}$ is zero on $\g a_2^{reg}$, one deduces from the Weyl integration formula that  there exist positive constants $C_\e$ (only depending of the choice of measures), such that  one has

 $$ \int_{\g q} F^+_{A,B}(X) f(X) dX=\sum_{\e\in\{ (++),(+-),(--)\}}C_\e \int_{\R^2} F^+_{A,B} (X_\e(u_1,u_2)) $$
$$ \times  \cM(f)(X_\e(u_1,u_2)) | u_1 u_2(\e_1 u_1^2-\e_2 u_2^2)| du_1 du_2.$$
 
By definition of $F^+_{A,B}$, there exist positive constants $C, C_1$ and $C_2 $ such that, for all $X_\e(u_1,u_2)\in  \O^{reg}$, one has 
$$ |(\e_1 u_1^2-\e_2 u_2^2) F^+_{A,B}(X_\e(u_1,u_2))|\leq C (C_1+|Ê\log\vert u_1\vert |)  (C_2+|Ê\log\vert u_2\vert |).$$

One deduces easily the corollary    from the following Lemma. \qed
 
\begin{lemma} Let $f\in\cD(\g q)$. Then there exist positive contants $C', C'_1, C'_2$ such that, for all $X_\e(u_1, u_2)\in \g q^{reg}$ one has
$$\vert \cM(f)(X_{ \e}(u_1, u_2))\vert \leq C' (C'_1+\big|\log \vert u_1\vert\big|)(C'_2+\big|\log \vert u_2\vert\big|).$$
\end{lemma}
\no{\it Proof.} Let $H=KNA$ be the Iwasawa decomposition of $H$ with $K=O(2)\times O(2)$, $N=N_0\times N_0$ where $N_0$ consists of 2 by 2 unipotent upper  triangular matrices and $A$ is the set of diagonal matrices in $H$. It is easy to see that the centralizer of $X$ in $H$ is the set of diagonal matrices $diag((\a,\b,\a,\b)$ with   $(\a,\b)\in(\R^*)^2$. Hence $H/Z_H(X)$ is isomorphic to $K\times  N \times \{diag(e^x,e^y,1,1); x,y\in \R\}$.

For $\xi  \in\R$, we set $n_\xi=\left(\begin{array}{cc} 1& \xi\\ 0 &1\end{array}\right)$.
We define the function $\tilde{f}$ by $\tilde{f}(X)=\int_K f(k\cdot X) dk$. Then, one has
  $$\cM(f)(X_{ \e}(u_1,u_2))=\vert \e_1 u_1^2-\e_2 u_2^2\vert \int_{\R^2}\big(\int_{\R^2} \tilde{f}(Y(u,\e,x,y,\xi,\eta)) d\xi d\eta  \big)dx dy$$
with $$Y(u,\e,x,y,\xi,\eta)=\Big( \left(\begin{array}{cc} n_\xi& 0\\ 0 &n_\eta\end{array}\right)diag(e^x, e^y,1,1)\Big)\cdot X_{\e,u}.$$
  Writing $Y(u,\e,x,y,\xi,\eta)=\left(\begin{array}{cc}0&Y\\ Z& 0 \end{array}\right)$, one has $$Y=\left(\begin{array}{cc}  u_1 e^x& -\eta  u_1 e^x+e^y\xi u_2\\ 0 & u_2 e^y \end{array}\right)\textrm{ and }Z=\left(\begin{array}{cc} \e_1u_1 e^{-x}& -\xi \e_1 u_1 e^{-x}+\eta \e_2 u_2e^{-y} \\ 0 & \e_2 u_2 e^{-y}\end{array}\right).$$\me
  
Since $f\in\cD(\g q)$, the function  $\tilde{f}$ has compact support in $\g q$. Identify $\g q$ with $\R^8$, there exists $T>0$ such that $\textrm{supp}(\tilde{f})\subset [-T,T]^8$. If  $\tilde{f}(Y(u,\e,x,y,\xi,\eta))\neq 0$ then we have  the following inequalities:
\begin{enumerate}\item $\vert u_1 e^{\pm x}\vert\leq T\quad\textrm{and }\quad\vert u_2 e^{\pm y}\vert\leq T$,
\item $\vert  -\eta  u_1 e^x+e^y\xi u_2\vert\leq T$,
 \item $\vert  -\xi \e_1 u_1 e^{-x}+\eta \e_2 u_2e^{-y}\vert\leq T$.
\end{enumerate}
Changing the variables $(\xi,\eta)$ in $(r,s)=( \xi u_2e^y-\eta  u_1 e^x,  -\xi \e_1 u_1 e^{-x}+\eta \e_2 u_2e^{-y})$, we obtain
the result.\qed

\no\underline{Remark.} By (\cite{HJ} Corollary 5.3.1), the function $F_{ana}$ defines an invariant eigendistribution on $\g q$. At this stage, we don't know if it is the case for the functions $F_{sing}$ and $F^+_{A, B}$. Indeed, the proof of Theorem \ref{resultatHJ} of \cite{HJ} is based on integration by parts using estimates of orbital integrals and  some of their derivates near non-zero semisimple elements of $\g q$. To determine if  $F_{sing}$ and $F^+_{A, B}$ are eingendistributions using the same method, we have to know the behavior of   derivates of  orbital integrals near $0$.

\end{document}